\documentclass{amsart}

\newtheorem{theorem}{Theorem}
\newtheorem{lemma}{Lemma}
\newtheorem{corollary}{Corollary}
\newtheorem{proposition}{Proposition}
\newcommand{\Z}{{\mathbf Z}}
\newcommand{\R}{{\mathbf R}}

\newcommand{\vol}{\operatorname{vol}}
\newcommand{\disc}{\operatorname{disc}}
\newcommand{\supp}{\operatorname{supp}}
\newcommand{\rad}{\operatorname{rad}}
\newcommand{\codim}{\operatorname{codim}}
\newcommand{\A}{\epsilon}
\newcommand{\K}{{\mathcal C}}
\newcommand{\LL}{{\mathcal L}}

\newcommand{\Cbin}{ r(r-1)/2 }

\newcommand{\Cmain }{ K}

\newcommand{\Ceps }{ K_1}
\newcommand{\Cdelta }{ K_2}
\newcommand{\CexpOne }{ K_3 }
\newcommand{\CchopQ }{ K_4}
\newcommand{\CepsProd }{K_5 }
\newcommand{\CexpTwo }{ K_3+1 }
\newcommand{\CchopC }{ K_6}
\newcommand{\CcompleteAgain }{ K_7}

\newcommand{\alphabound}{\mbox{$\lceil 1/2 +
    \frac{ \sqrt{\omega(q)} } {7 \log_2 p} \rceil $}}

\begin{document}
\author{P\"ar Kurlberg}
\title{The distribution of spacings between quadratic residues, II}
\begin{abstract}
  We study the distribution of spacings between squares in $\Z/Q\Z$ as
  the number of {\em prime} divisors of $Q$ tends to infinity. In
  \cite{part1} Kurlberg and Rudnick proved that the spacing
  distribution for {\em square free} $Q$ is Poissonian, this paper
  extends the result to {\em arbitrary} $Q$.
\end{abstract}
\address{Raymond and Beverly Sackler School of Mathematical Sciences,
Tel Aviv University, Tel Aviv 69978, Israel} \date{August 4, 1998}
\thanks{Supported in part by grants from the Israel Science Foundation
and by the EC TMR network "Algebraic Lie Representations", EC-contract
no ERB FMRX-CT97-010}

\maketitle

\section{Introduction} 

This paper studies the distribution of spacings between squares in
$\Z/Q\Z$ as $\omega(Q)$, the number of {\em prime} divisors of $Q$,
tends to infinity. In \cite{part1} Kurlberg and Rudnick proved that
the spacing distribution for square free $Q$ is {\em Poissonian},
i.e., the same as for a sequence of independent uniformly distributed
real numbers in the unit interval. The purpose of this paper is to
extend the result to arbitrary $Q$.

The spacing distribution is defined as follows: Let $X_Q \subset \{ 0,
1, \ldots, Q-1 \}$ be a set of representatives of the squares in
$\Z/Q\Z$. Order the $N_Q$ elements of $X_Q$ so that $x_1 < x_2 <
\ldots < x_{N_Q}$ and form the {\em normalized consecutive spacings}
$y_i = (x_{i+1}-x_i)/s$ where $s=(x_{N_Q}-x_1)/N_Q$ is the mean
spacing. By putting point mass $(N_Q-1)^{-1}$ at each $y_i$ we obtain
a probability distribution with mean one, and we can now study the
limiting distribution as $\omega(Q) \rightarrow \infty$.

For prime $Q \rightarrow \infty$ the mean spacing is constant and
Davenport \cite{dav} has proved that the normalized spacing
distribution is a sum of point masses at half integers $k/2$ with
weight $2^{-k}$, and it is easy to see that the same holds true for
prime powers. In the highly composite case the mean spacing tends to
infinity since $s$ roughly equals $2^{\omega(Q)}$. Hence there is a
chance that the limiting distribution has continuous support.
Davenport's result in a sense suggests that quadratic residues behave,
at least with respect to spacing statistics, like independent fair
coin flips. This together with the heuristic that primes are
independent suggests that the limiting distribution for highly
composite $Q$ should be Poissonian, i.e., the probability
density function of the (normalized) spacing to the next square should
be given by $P(s) = e^{-s}$.

The definition of the level spacing distribution involves ordering the
elements in $X_Q$. In terms of analysis, ordering is a complicated
operation, and it is not so easy to study the level spacings directly.
However, using a combinatorial argument one can recover the level
spacings from the knowledge of all {\em $r$-level correlations.} (For
instance, see lemma~14 of \cite{part1}.)

Fix an integer $r \geq 2$. The {\em $r$-level correlation} is defined
as follows: let $\K \subset \R^{k-1}$ be a convex set such that
$(x_1-x_2, x_2-x_3, \ldots, x_{k-1}-x_{x}) \in \K$ implies that $x_i
\neq x_j$ for $i \neq j$. The reason for this condition is that we
want to avoid the self correlation of a point with itself. The
$r$-level correlation with respect to $\K$ is given by:
$$
R_r( \K, Q) = 
\frac{1}{N_Q} 
\sum_{h \in s \K \cap \Z^{r-1}}
N(h,Q)
$$
where $N(h,Q)$ is the number of solutions in squares $s_1,
\ldots ,s_r$ of the equations 
$$
s_{i+1}-s_i \equiv h_i \mod Q, \ \ i=1, \ldots, r-1.
$$

The main result of this paper is the following:
\begin{theorem} \label{t:main} With $\K$ as above there exists a
  constant $\Cmain>0$, depending only on $r$ and $\K$, such that
  $$
  R_r( \K, Q) = \vol( \K ) + O(\exp( -\Cmain \sqrt{\omega(Q)} )).
  $$
\end{theorem}
As is well known (for instance, see lemma~14 of \cite{part1}) this
implies:  
\begin{corollary}
  The limiting spacing distribution of squares in $\Z/Q\Z$ as
  $\omega(Q) \rightarrow \infty$ is Poissonian.
\end{corollary}

{\em Remark:} In the special case that the exponents of the primes
dividing $Q$ are bounded then the methods in \cite{part1} can be
generalized. For the general case one can try to truncate $Q$, i.e.
replace $Q$ by $\tilde{Q}$ in such a way that that the growth of the
exponents is controlled.  However, new ideas are needed in order to
justify that the errors introduced by truncating $Q$ and cutting off
certain divisor sums are not too big. Because of the ensuing
complications the bound on the error term in theorem~\ref{t:main} is
only of sub-exponential decay in $\omega(Q)$, whereas in theorem~1 of
\cite{part1} the bound decays exponentially.

{\em Contents of the paper:} In section 2 we set up the necessary
notation, and in section 3 we show that the decomposition of $N(h,p)$
used in \cite{part1} is valid for prime powers. Squares that are
distinct modulo $Q$ are not necessarily distinct modulo $p$, and in
section 4 we briefly recall some properties of this modulo $p$
degeneracy  and its relation to lattices and M\"oebius inversion.
Section 5 deals with truncating $Q$, i.e., lowering the exponents of
the primes dividing $Q$, as well as truncating sums over sets of
lattices and divisors of $Q$.  In section 6 we use the previous
results to show that a periodicity heuristic is valid, using it we
prove theorem~\ref{t:main}. Finally, in the appendix we collect some
lemmas on divisor sums used throughout the text.

\section{Notation}

For $n$ an integer we let $\omega(n)$ be the number of {\em prime}
divisors of $n$. When writing $p|n$ we will {\em always} refer to a
prime divisor of $n$. Let $Q = \prod_{p|q} p^{\alpha_p}$ where $q =
\rad(Q)$ is the square free part of $Q$. (Note that
$\omega(Q)=\omega(q)$.) Put $\tilde{Q}= \prod_{p|q}
p^{\tilde{\alpha}_p}$ where $\tilde{\alpha}_p \leq \alpha_p$ are to be
picked later. If $c|q$ we let $C = \prod_{p|c} p^{\alpha_p}$ and
$\tilde{C} =\prod_{p|C} p^{\tilde{\alpha}_p}$.

In what follows we will use the following convention: If a function,
say $f$, is defined for prime arguments we let
$$f(c) = \prod_{p|c} f(p).$$
If the function is defined for prime powers, let 
$$f(C) = \prod_{p|c} f(p^{\alpha_p})$$
and
$$f(\tilde{C}) = \prod_{p|c} f(p^{\tilde{\alpha}_p}).$$
For instance, we let $\sigma(p) = 1+p^{-1}$; by the above convention
$\sigma(q) = \prod_{p|q} \sigma(p) = \sum_{c|q} c^{-1}$.

We let $s=Q/N_Q$ denote the mean spacing. (This is slightly different
from what is used in the introduction, but in the limit $\omega(Q) \rightarrow
\infty$ the two definitions agree.) It is easy to see that
$N_{p^k}=p^k \frac{\sigma(p)}{2} \left( 1+O(p^{-2}) \right)$, with the
error term always positive, and therefore
$$
2^{\omega(q)(1-\epsilon)} \ll
\frac{2^{\omega(q)}}{\sigma(q)} 
\leq s \ll
\frac{2^{\omega(q)}}{\sigma(q)}
\leq 2^{\omega(q)}.
$$
Finally, we let 
$$F(q,t) =\sum_{p|q} p^{-t}.$$

\section{Analyzing $N(h,Q)$}

Since $x$ is a square modulo $Q$ if and only if it is a square
modulo $P$ for all $P|Q$, we see that $N(h,Q)$ is multiplicative. For
primes we have:
\begin{lemma}
  \label{l:N-prime}
  We can write 
  $$
  N(h,p) = \frac{\Delta(h,p)}{2^r} p \left( 1 + \A(h,p) \right)
  $$
  where $$\A(h,p) \ll_r  p^{-1/2} \qquad \mbox{as }  p \rightarrow
  \infty,$$ 
  and $\Delta(h,p) = 2^k$ for some $0 \leq k \leq r$. 
\end{lemma}
\begin{proof} See proposition~4 in \cite{part1}.
\end{proof}
{\em Remark:} For $p \neq 2$ the bound on $\A(h,p)$ follows from the Weil
bounds on the number of points on curves over finite fields. For
$p=2$ the curve is highly singular, but the bound holds trivially by
choosing a large enough constant. 

For prime powers Hensel's lemma can be used to lift
solutions. However, there are 
complications due to singularities arising for certain choices of
$h$. (Note that if all points were smooth, then $N(h,p^k) =
p^{k-1}N(h,p)$.) The following lemma shows that there are few
solutions that do not lift. 
\begin{proposition} 
  \label{p:N-chop-bound}
  If $b \geq a$ then
  $$ | N(h, p^b)- p^{b-a} N(h,p^a)| \ll_r p^{b-a}.$$
\end{proposition} 
\begin{proof}
Recall that $N(h,p^b)$ is the number of solutions in squares $s_1,
\ldots ,s_r$ of the equations 
$$
s_{i+1}-s_i \equiv h_i \mod p^b, \ \ i=1, \ldots, r-1,
$$
which we may rewrite as
$$
y_i^2 \equiv t_i + x^2 \mod p^b, \ \ i=1, \ldots, r-1,
$$
where $t_i = \sum_{j=1}^i h_j$ and we think of $x$ as a preferred
parameter. For most values of 
$x$, the equations in $y_i$ are {\em smooth}, and Hensel's
lemma can be applied to lift solutions modulo $p$ to
solutions for arbitrary high powers. However, at the non-smooth points
the analysis is more involved. 

Assume first that $p \neq 2$. For the pair correlation we get the
equation  
$$
y^2 \equiv x^2+t \mod p^b.
$$ 
If $x^2+t \not \equiv 0 \mod p$, we're in the smooth case. If not,
then we can write $ x^2+t = u p^k$, where $u$ is invertible modulo $p$
and $1 \leq k \leq b$. Now, $y^2 \equiv u p^k \mod p^b$ has a solution
iff $k$ is even and $u$ is a square in $F_p$, or $k=b$. Thus, if $y^2
\equiv x^2+t \mod 
p^a$ has a solution which does not lift, this implies that $k \geq a$.
The $x^2$ for which solutions cannot be lifted are contained in the
$(p^b)$-cosets generated by $-t + (p^a)$, and there are at most
$|(p^a)/(p^b)|=p^{b-a}$ such elements in $\Z/p^b\Z$. For $r \geq 3$ we
observe that the ``bad'' $x^2$ are contained in the $p^b$-cosets
generated by $\cup_{i=1}^{r-1} \left( -t_i + (p^a) \right)$, and there
are at most $(r-1) p^{b-a}$ such $x^2$.

For $p=2$ the difference is that a unit has to be a square modulo
$8$ in order to be a dyadic square, and we therefore lose a
factor of $4$ when bounding the number of ``bad'' squares.
\end{proof}

\begin{corollary}
  \label{p:N-prime-power}
  We can write 
  $$
  N(h,p^k) = \frac{\Delta(h,p)}{2^r} p^k \left( 1 + \A(h,p^k) \right)
  $$
  where $$\A(h,p^k) \ll_r  p^{-1/2}.$$
\end{corollary}

\begin{corollary}
\label{c:eps-bound}
There exists $\Ceps>0$ such that 
$$
\A(h,C) \leq \Ceps^{\omega(c)} c^{-1/2}
$$
for all $c|q$.
\end{corollary}

\section{$\Delta$, lattices and M\"obius inversion}
In this section we briefly explain $\Delta(h,p)$, which measures
how many extra solutions in squares $s_i$ of the system 
$$
s_{i+1}-s_i \equiv h_i \mod p, \ \ i=1, \ldots, r-1,
$$
there are. (For full details see section~4.1 in \cite{part1}.) For the
pair correlation 
it works as follows: if $h \not \equiv 0 \mod
p$ then there are roughly 
$N_p/2 \simeq p/4$ solutions of $s_2 \equiv s_1 + h \mod p$ since the
``probability'' of $s_1 + h$ being a square modulo $p$ is roughly
$1/2$. However, if $h \equiv 0 \mod p$ there is degeneracy;
$s_2=s_1$ is {\em automatically} a square. Hence there are $N_p \simeq
p/2$ solutions in this case, and $\Delta(0,p)=2$ is the corresponding 
correction factor. 

For the higher correlations the ``probability'' of $s_1, s_1+h_1,
s_1+h_1+h_2, \ldots, s_1+h_1+ \ldots h_{r-1}$ all being squares is
roughly $2^{-r}$, {\em unless} the values of the $h_i$'s forces some
of the $s_i$'s to be equal. More precisely, if we let ${\mathcal H}_p$
be the union of linear subspaces in $(\Z/p\Z)^{r-1}$ that corresponds
to some $s_i$'s being equal, then the condition for degeneracy translates
into $h$ lying in a unique smallest linear subspace $H \in {\mathcal
  H}_p$, and the corresponding correction factor is
$2^{\codim(H)}$.

Using M\"obius inversion we can express the function $\Delta(h,p) =
2^{\codim(H)}$ as a linear combination of characteristic functions of
the linear subspaces in ${\mathcal H}_p$. Pulling the subspaces in
$(\Z/p\Z)^{r-1}$ back to $\Z^{r-1}$ gives a set of lattices $\LL_p$,
and $\Delta(h,p) = \sum_{L \in \LL_p} \lambda(L) \delta_L(h)$ where
$\delta_L$ is the characteristic function of the lattice and
$\lambda(L)$ are certain coefficients (see section 4.1 in
\cite{part1}.) Note that the set of values $\{ \lambda(L) \ | \ L \in
\LL_p \}$ is independent of $p$, and that $\lambda(L)=1$ for the the
maximal lattice $L=\Z^{r-1}$.

For divisors $c|q$, we then have  
$$
\Delta(h,c) = \prod_{p|c} \Delta(h,p) =
\prod_{p|c} \left(
  \sum_{L_p \in \LL_p}
  \lambda(L_p) \delta_{L_p}(h)
\right)
=
\sum_{g|c} \sum_{ \substack{ L \in \LL \\ \supp(L)=g  }  }
\lambda(L) \delta_L(h)
$$ where the inner sum is over the collection $\LL$ of all lattices of
the form $\cap_{p|q} L_p, L_p \in \LL_p$, the coefficient $\lambda
(\cap_{p|q} L_p) $ is given by $\prod_{p|q} \lambda(L_p)$, and where
we let $\supp(L)$, the support of $L$, be the square free part of the
discriminant $\disc(L)$. (The discriminant is as usual the volume of
the fundamental domain of the lattice.)

{\em Remark:} If $L = \cap_{p|q} L_p$ and $p$ does not divide
$\supp(L)$ then $L_p = \Z^{r-1}$ and $\lambda(L_p)=1$. Consequently,
$L = \cap_{p|\supp(L)} L_p$ and $\lambda(L) =\prod_{p|\supp(L)}
\lambda(L_p)$.

For future reference we have the following lemmas:
\begin{lemma} 
  \label{l:lattice-bound}
  $$
  |\lambda(L)|  \ll_{r} \supp(L)^{\epsilon},
  $$  
  and
  $$
  \sum_{ \substack{ L \in \LL \\ \supp(L)=g  }  } 1
  \ll_{r} g^{\epsilon}.
  $$
\end{lemma}
\begin{proof}
Immediate by the previous remark and the fact that
$\lambda(L_p),|\LL_p| \ll_r 1$.  
\end{proof}

\begin{lemma}
\label{l:delta-sup} The following bound holds:
$$
\Delta(h,c) \ll c^{\epsilon}.
$$
\end{lemma}
\begin{proof}
  By the previous remark there exists a constant $\Cdelta$, depending
  only on $r$, such that $\Delta(h,p) \leq \Cdelta$. Hence
  $\Delta(h,c) = \prod_{p|c} \Delta(h,p) \leq \Cdelta^{\omega(c)} \ll
  c^{\epsilon}$.
\end{proof}

By assumption the convex set $\K$ has empty intersection with the
linear subspaces, or walls, corresponding to $ \sum_{ i \leq j \leq k
  } h_j = 0$ for $1\leq i \leq k \leq r-1$. The lattices in $\LL$
correspond to integer points that are congruent to the walls modulo
some divisor of $Q$. Thus, if the support of a lattice is sufficiently
large compared to the size of $s \K$ we expect it to have empty
intersection with $s \K$, and this is in fact true:

\begin{lemma} 
  \label{l:empty-intersection}
  If $\K$ does not intersect with the walls and $\supp(L) \gg_{\K}
  s^{\Cbin}$ then $s\K \cap L = \emptyset$.  
\end{lemma}
\begin{proof}
See lemma~7 in \cite{part1}. 
\end{proof}

If $L \subset \R^n$ is a lattice and $X \subset \R^n$ is a set with
nice boundary, for instance if $X$ is convex, then it is well
known that the number of lattice points in $t \cdot X $ equals $t^n
\frac{\vol(X)}{\disc(L)} + O_{X,L}(t^{n-1})$, where the error term depends
on the set $X$ {\em and} the lattice $L$.  The Lipschitz principle
(Davenport \cite{dav2}, Schmidt \cite{schmidt}) allows us to bound the
error uniformly with respect to {\em integer lattices} $L \subset
\Z^n$:
\begin{proposition}
\label{p:lipschitz} 
Let $L \subset \Z^n$ be a lattice of discriminant $\disc(L)$, and $\K$
a convex set. Suppose that $\K$ lies in a ball of radius $R$. Then
$$
\#(L \cap \K) = \frac{\vol(\K)}{\disc(L)} + O_{\K}(R^{n-1})
$$
where the error term only depends on $\K$.
\end{proposition}
\begin{proof}
For details see lemma~16 in \cite{part1}.
\end{proof}

The following bounds the sum of $\Delta(h,q)$ over all integer points in $s
\K$.  
\begin{lemma}
\label{l:Delta-sum}
  $$
  \sum_{ h \in s \K \cap \Z^{r-1}} \Delta(h,q)
  \ll_r
  s^{r-1} \exp \left( O \left( \log \log \left(\omega(q) \right)
    \right) \right)
  $$
\end{lemma}
\begin{proof}
Rewriting the sum using M\"obius inversion and using
proposition~\ref{p:lipschitz} we get 
$$
\sum_{ h \in s \K \cap \Z^{r-1}} \Delta(h,q)
=
\sum_{g|q}
\sum_{ \substack{ L \in \LL \\ \supp(L)=g } }
\lambda(L) 
\sum_{ h \in s \K \cap L} 1
$$

$$ =
\sum_{g|q}
\sum_{ \substack{ L \in \LL \\ \supp(L)=g } }
\lambda(L) 
\left(
  \frac{ \vol(s \K) }{ \disc(L) } +
  O(s^{r-2})
\right).  
$$
By lemma~\ref{l:empty-intersection} we may assume that $g \leq
s^{\Cbin}$, and we may estimate the terms involving $O(s^{r-2})$ by 
$$
\sum_{ \substack{  g|q  \\ g \leq s^{\Cbin}    } }
\sum_{ \substack{ L \in \LL \\ \supp(L)=g } }
|\lambda(L)|   s^{r-2}
\ll
s^{r-2} 
\sum_{ \substack{  g|q  \\ g \leq s^{\Cbin}    } }
g^{2 \epsilon}
$$
using  lemma~\ref{l:lattice-bound}. But 
$$
\sum_{ \substack{  g|q  \\ g \leq s^{\Cbin}    } }
g^{2 \epsilon}
\ll
s^{ 2 \epsilon \Cbin} 
\sum_{ \substack{  g|q  \\ g \leq s^{\Cbin}    } }
1
\ll 
s^{\epsilon'}
$$
by lemma~\ref{l:div-sum-bound} and hence the error terms only
contribute $O(s^{r-2+\epsilon})$. The main term 
$$
\sum_{ \substack{  g|q  \\ g \leq s^{\Cbin}    } }
\sum_{ \substack{ L \in \LL \\ \supp(L)=g } }
\lambda(L) 
\frac{ \vol(s \K) }{ \disc(L) } 
$$
is trivially bounded by 
$$
\vol(s \K) \sum_{ g|q } 
\sum_{ \substack{ L \in \LL \\ \supp(L)=g } }
\frac{ |\lambda(L)|}
{ \disc(L) }.
$$

Using multiplicativity once more we get 
$$
\sum_{  g|q   }
\sum_{ \substack{ L \in \LL \\ \supp(L)=g } }
\frac{ |\lambda(L)|  }{ \disc(L) } 
=
\prod_{  p|q   } \left(
  \sum_{  L \in \LL_p }
  \frac{ |\lambda(L)|  }{ \disc(L) }
\right).
$$ Now, if $L \in \LL_p$ then $\disc(L)$ is a power of $p$, and unless
$L = \Z^{r-1}$ the power is $\geq 1$. Since $\lambda(\Z^{r-1})=1$ we
see that $\sum_{ L \in \LL_p } \frac{ |\lambda(L)| }{ \disc(L) } = 1 +
O(p^{-1})$. Recalling that the number of lattices in $\LL_p$ and the
set of values $\{ \lambda(L) \ | \ L \in \LL_p \}$ are independent of
$p$ we see that the error is uniform in $p$. Thus
$$
\sum_{  g|q   }
\sum_{ \substack{ L \in \LL \\ \supp(L)=g } }
\frac{ |\lambda(L)|  }{ \disc(L) } 
\ll
\prod_{  p|q   } \left(  1 + O(p^{-1}) \right) 
\ll
\exp \left( \sum_{  p|q   } O(p^{-1}) \right) 
=
\exp \left( O \left( F_q(1) \right)  \right)
$$
By lemma~\ref{l:f-sub-q-bound}, $F_q(1) = O \left( \log \log
  \left( \omega(q) \right) \right) $ and we are done. 
\end{proof}

\section{Truncations}
In order to use periodicity in section 6 we will need to control
the error when we replace $Q = \prod_{p|q} p^{\alpha_p}$ by $\tilde{Q}
= \prod_{p|q} p^{\tilde{\alpha}_p}$, where
$$
\tilde{\alpha}_p = \min( \alphabound, \alpha_p).
$$ 
We will also need to show that sums over large divisors and
lattices are small.  

\subsection{Truncating $Q$}
First note that if $\tilde{\alpha}_p < \alpha_p$ then
$p^{-\tilde{\alpha}_p} \leq p^{-1/2} \exp \left( -\frac{\sqrt{
      \omega(q) }}{\CexpOne} \right)$ for $\CexpOne > \frac{7}{\log
  2}$.  The following shows that we are not committing too large of an
error when we truncate $Q$.
\begin{proposition}
  \label{p:main-chop}
  There exists $\CchopQ>0$ such that 
  $$
  R_r( \K, Q) = 
  s \sum_{h \in s \K \cap \Z^{r-1}}
  \frac{N(h,\tilde{Q})}{\tilde{Q}}
  + O\left( \exp( -\CchopQ \sqrt{\omega(q)} ) \right). 
  $$
\end{proposition}
\begin{proof}

First we prove the following claim:
$$
| \frac{N(h,Q)}{Q} - \frac{N(h,\tilde{Q})}{\tilde{Q}}|
=
\frac{1}{Q} 
\left| 
  \prod_{p|q} N(h, p^{\alpha_p}) -
  \prod_{p|q} N(h, p^{\tilde{\alpha}_p}) p^{\alpha_p - \tilde{\alpha}_p}
\right|
$$
$$
\ll 
\frac{\Delta(h,q)}{2^{r\omega(q)}} 
\exp \left( \CepsProd \sqrt{\frac{\omega(q)}{\log
    \omega(q)}}   -  \frac{\sqrt{\omega(q)}}{\CexpTwo} \right). 
$$
Letting $A_p = N(h, p^{\alpha_p})$ and $B_p = N(h, p^{\tilde{\alpha}_p})
p^{\alpha_p - \tilde{\alpha}_p}$ we have $|A_p-B_p| \ll_r p^{\alpha_p
  - \tilde{\alpha}_p}$ by proposition~\ref{p:N-chop-bound}. We may assume
that $B_p$ is nonzero for all $p$ since $B_p =0$ implies
that $A_p =0$ (there are no solutions to lift), and if $A_p
=B_p=0$ the bound holds trivially. Now,  
$$
|\prod_{p|q} A_p - \prod_{p|q} B_p| =
\left( \prod_{p|q} B_p \right)
\left|
  \prod_{p|q} \left(
    1 + \frac{A_p-B_p}{B_p}
  \right)
  -1
\right|
$$
$$
=
\left( \prod_{p|q} B_p \right) \left| 
  \exp \left( \sum_{p|q}
    \log \left(  \frac{A_p-B_p}{B_p} +1  \right)
  \right)
  -1
\right|
\ll
\left( \prod_{p|q} B_p \right) 
\sum_{p|q} \left| \frac{A_p-B_p}{B_p} \right|
$$
Thus
$$
| \frac{N(h,Q)}{Q} - \frac{N(h,\tilde{Q})}{\tilde{Q}}|
\ll
\sum_{  p_0 | q         } 
\left(
  \prod_{ p| \frac{q}{p_0} }
  \frac{B_p}{p^{\alpha_p}}
\right)
\frac{|A_{p_0}-B_{p_0}|}{p_0^{\alpha_{p_0}}}
$$
$$
\ll
\sum_{  p_0 | q         } 
\left(
  \prod_{ p| \frac{q}{p_0} }
  \frac{N(h,p^{\tilde{\alpha}_p})}{p^{\tilde{\alpha}_p}}
\right)
p_{0}^{-\tilde{\alpha}_{p_0}}
$$
$$
\ll
\sum_{  p_0 | q         } 
\prod_{ p| \frac{q}{p_0} }
\frac{\Delta(h,p)}{2^r} (1+ \Ceps p^{-1/2})
\exp \left( - \frac{\sqrt{\omega(q)}}{\CexpOne} \right)
$$
by corollary~\ref{c:eps-bound} and 
since we can assume that $\tilde{\alpha}_p < \alpha_p$. (If they are
equal then $A_p=B_p$.) This is in turn bounded by
$$
\omega(q) 
\exp \left( - \frac{\sqrt{\omega(q)}}{\CexpOne} \right)
\prod_{ p|q} (1+ \Ceps p^{-1/2})
$$
since $\frac{\Delta(h,p)}{2^r} \leq 1$. By lemma~\ref{l:chop-complete} 
$$
\prod_{ p|q} (1+ \Ceps p^{-1/2}) \ll 
\exp \left( \CepsProd \sqrt{  \frac{ \omega(q)}{\log \omega(q)}} \right)
$$
and we have proved the claim.

Summing over all $h$ and applying lemma~\ref{l:Delta-sum}
gives that  
$$
s \sum_{h \in s \K \cap \Z^{r-1}}
\frac{\Delta(h,q)}{2^{r\omega(q)}}
\ll
\frac{s^r}{2^{r\omega(q)}} 
\exp \left( O(\log \log \omega(q)) \right) 
\ll
\exp \left( O(\log \log \omega(q)) \right).
$$
Hence 
$$
s \sum_{h \in s \K \cap \Z^{r-1}}
| \frac{N(h,Q)}{Q} - \frac{N(h,\tilde{Q})}{\tilde{Q}}|
$$
$$
\ll
\exp \left( O \left( \log\log(\omega(q)) \right) +  \CepsProd
\frac{\sqrt{\omega(q)}}{\log \omega(q)}   - \frac{\sqrt{\omega(q)}}{\CexpTwo}
\right).   
$$
and we are done as $s=Q/N_Q$ and thus
$$
R_r( \K, Q) = \frac{1}{N_Q} \sum_{h \in s \K \cap \Z^{r-1}} N(h,Q) =
s \sum_{h \in s \K \cap \Z^{r-1}} \frac{N(h,Q)}{Q}. 
$$
\end{proof}

\begin{corollary}
  There exists $\CchopQ>0$ such that 
  $$
  R_r( \K, Q) = 
  \frac{s}{ 2^{r\omega(q)}} 
  \sum_{c|q}
  \sum_{h \in s \K \cap \Z^{r-1}}
  \Delta(h,q) 
  \A(h,\tilde{C})
  + O\left( \exp( -\CchopQ \sqrt{\omega(q)} ) \right). 
  $$
\end{corollary}
\begin{proof}
Immediate since 
$$
\frac{N(h,\tilde{Q})}{\tilde{Q}} =  
\prod_{p|q}
\frac{\Delta(h,p)}{2^r} \left( 1 + \A(h,p^{\tilde{\alpha}_p})
\right) 
=
\frac{\Delta(h,q)}{2^{r\omega(q)}} \sum_{c|q} \A(h,\tilde{C}).
$$
\end{proof}

The choice of $\tilde{\alpha}_p$ also gives some control
of the size 
of $\tilde{C}$ when $\omega(c) \leq \sqrt{\omega(q)}$:
\begin{lemma}
  \label{l:bound-d-tilde}
  If $c|q$ and $\omega(c) \leq \sqrt{\omega(q)}$ then $\tilde{C}
  \leq c^{3/2} s^{1/6}$. 
\end{lemma}
\begin{proof} We have $\tilde{\alpha}_p \leq 3/2 + \frac{ \sqrt{\omega(q)} }
  {7 \log_2 p}  $ and thus
$$
\log_2 \tilde{C} \leq 
3/2 \log_2 c +
\sum_{p|c}  \log_2 p \frac{ \sqrt{\omega(q)} } {7 \log_2 p} 
\leq
3/2 \log_2 c +
\omega(c) \frac{ \sqrt{\omega(q)} } {7} 
$$
$$ \leq 3/2 \log_2 c +\frac{ \omega(q) } {7} .
$$
Exponentiating we get $\tilde{C} \leq c^{3/2} 2^{\omega(q)/7}
\ll c^{3/2} s^{1/7+\epsilon} \ll 
c^{3/2} s^{1/6}$.
\end{proof}

\subsection{Truncating divisor sums}
In order to use periodicity in section 6 we need the product of
$\tilde{C}$ 
and certain discriminants of lattices to be small. We will prove that
the contribution of terms where this is not the case is
negligible. First we show that $c$ with many divisors, or of large
size, can be neglected.
\begin{lemma} \label{l:chop-c} There exists $\CchopC>0$ such that  
$$
\frac{s}{ 2^{r\omega(q)}} 
\sum_{c|q}
\sum_{h \in s \K \cap \Z^{r-1}}
\Delta(h,q) 
\A(h,\tilde{C})
$$
$$
=
\frac{s}{ 2^{r\omega(q)}} 
\sum_{ \substack{ c|q \\ c \leq s^{1/3} \\ \omega(c) \leq \sqrt{\omega(q)} } }
\sum_{h \in s \K \cap \Z^{r-1}}
\Delta(h,q) 
\A(h,\tilde{C})
+ O \left( \exp \left( -\frac{\CchopC}{2} \sqrt{\omega(q)}  \right) \right). 
$$
\end{lemma}

\begin{proof}
By corollary~\ref{c:eps-bound},
$$
\sum_{h \in s \K \cap \Z^{r-1}}
\left| 
  \Delta(h,q) \A(h,\tilde{C})
\right|
$$
$$
\ll
c^{-1/2} \Ceps^{\omega(c)} 
\sum_{h \in s \K \cap \Z^{r-1}}
\Delta(h,q). 
$$
Moreover, 
$$
\frac{s}{ 2^{r\omega(q)}} 
\sum_{h \in s \K \cap \Z^{r-1}}
\Delta(h,q)
\ll 
\frac{s^r}{2^{r\omega(q)}} \exp \left( O(\log \log \omega(q)) \right)
$$
$$
\ll
\exp \left( O(\log \log \omega(q)) \right)
$$
by lemma~\ref{l:Delta-sum} and the result now follows from
lemma~\ref{l:chop-complete}.
\end{proof}

\subsection{Truncating lattice sums}
Writing $\Delta(h,q)= \Delta(h,c) \cdot \Delta(h,\frac{q}{c})$ and expanding
the second term 
$$
\Delta(h,\frac{q}{c})=
\sum_{g | \frac{q}{c}}
\sum_{ \substack{ L \in \LL \\ \supp(L) = g } }
\lambda(L) \delta_L(h)
$$
we get 
$$
\frac{s}{ 2^{r\omega(q)}} 
\sum_{ \substack{ c|q \\ c \leq s^{1/3} \\ \omega(c) \leq \sqrt{\omega(q)} } }
\sum_{h \in s \K \cap \Z^{r-1}}
\Delta(h,q) 
\A(h,\tilde{C})
$$
$$
=
\frac{s}{ 2^{r\omega(q)}} 
\sum_{ \substack{ c|q \\ c \leq s^{1/3} \\ \omega(c) \leq \sqrt{\omega(q)} } }
\sum_{g | \frac{q}{c}}
\sum_{ \substack{ L \in \LL \\ \supp(L) = g } }
\lambda(L)
\sum_{h \in s \K \cap L}
\A(h,\tilde{C}) \Delta(h,c).
$$
We now show that lattices with large discriminants can be
neglected:
\begin{lemma} \label{l:chop-l}
$$
\frac{s}{ 2^{r\omega(q)}} 
\sum_{ \substack{ c|q \\ c \leq s^{1/3} \\ \omega(c) \leq \sqrt{\omega(q)} } }
\sum_{g | \frac{q}{c}}
\sum_{ \substack{ L \in \LL \\ \supp(L) = g } }
\lambda(L)
\sum_{h \in s \K \cap L}
\A(h,\tilde{C}) \Delta(h,c) 
$$
$$
=
\frac{s}{ 2^{r\omega(q)}} 
\sum_{ \substack{ c|q \\ c \leq s^{1/3} \\ \omega(c) \leq \sqrt{\omega(q)} } }
\sum_{ \substack{ g | \frac{q}{c}  \\ g \leq s^{\Cbin} } }
\sum_{ \substack{ L \in \LL \\ \supp(L) = g \\ \disc(L) \leq s^{1/3}  } }
\lambda(L)
\sum_{h \in s \K \cap L}
\A(h,\tilde{C}) \Delta(h,c) 
+
O(s^{-1/3 + \epsilon}).
$$
\end{lemma}
\begin{proof}
By lemma~\ref{l:empty-intersection}, we may assume that $g \leq
s^{\Cbin}$. By corollary~\ref{c:eps-bound}, and
lemma~\ref{l:delta-sup} $\A(h,\tilde{C}) \Delta(h,c) \ll c^{-1/2
  +\epsilon }$, thus 
$$
\sum_{h \in s \K \cap L}
\A(h,\tilde{C}) \Delta(h,c) 
\ll
c^{-1/2 +\epsilon } 
\sum_{h \in s \K \cap L}
1.
$$ 
By proposition~\ref{p:lipschitz} 
$$
\sum_{h \in s \K \cap L} 1 = 
\frac{\vol(s \K)}{\disc(L)}+O(s^{r-2}).$$ 
Now, 
$$
\frac{s}{ 2^{r\omega(q)}} 
\sum_{ \substack{ c|q \\ c \leq s^{1/3} \\ \omega(c) \leq \sqrt{\omega(q)} } }
c^{-1/2 +\epsilon } 
\sum_{ \substack{ g | \frac{q}{c}  \\ g \leq s^{\Cbin}    }}
\sum_{ \substack{ L \in \LL \\ \supp(L) = g \\ \disc(L) \geq s^{1/3}  } }
|\lambda(L)|
\left(
  \frac{\vol(s \K)}{\disc(L)}
  + O(s^{r-2})
\right)
$$
$$
\ll 
\frac{s}{ 2^{r\omega(q)}} 
\sum_{ \substack{ c|q \\ c \leq s^{1/3} \\ \omega(c) \leq \sqrt{\omega(q)} } }
c^{-1/2 +\epsilon } 
\sum_{ \substack{ g | \frac{q}{c}  \\ g \leq s^{\Cbin}        }}
g^{2\epsilon}
\left(
  \frac{\vol(s \K)}{s^{1/3}}
  + O( s^{r-2} )
\right)
$$
by lemma~\ref{l:lattice-bound}. Thus we need to show that 
$$
\frac{\vol(\K) s^r}{ 2^{r \omega(q)}} 
\sum_{ \substack{ c|q \\ c \leq s^{1/3} \\ \omega(c) \leq \sqrt{\omega(q)} } }
c^{-1/2 +\epsilon } 
\sum_{ \substack{ g | \frac{q}{c}  \\  g \leq s^{\Cbin}  }}
O(s^{-1/3})
= O(s^{-1/3+\epsilon}).
$$ 
By lemma~\ref{l:div-sum-bound}
$$
\sum_{ \substack{ g | \frac{q}{c}  \\  g \leq s^{\Cbin}  }} s^{-1/3}
\ll s^{-1/3+\epsilon},
$$
and 
$$
\sum_{ \substack{ c|q \\ c \leq s^{1/3} \\ \omega(c) \leq \sqrt{\omega(q)} } }
c^{-1/2 +\epsilon } 
\ll
\sum_{ \substack{ c|q \\ c \leq s^{1/3} \\ \omega(c) \leq
    \sqrt{\omega(q)} } } O(1)
\ll s^{\epsilon}.
$$ 
Since $\frac{\vol(\K) s^r}{ 2^{r \omega(q)}} \leq \frac{\vol(\K)}{
  \sigma(q)^{r}} \leq 1$ we see that the sum over $L$ such that
$\disc(L) \geq s^{1/3}$ is $O(s^{-1/3+\epsilon})$.
\end{proof}

\section{Periodicity}

We are now in the position of using periodicity of $\A(h,\tilde{C})
\Delta(h,c)$ modulo $\tilde{C}$, i.e. if $\disc(L) \cdot \tilde{C}
\leq s$ then
$$
\sum_{h \in s \K \cap L}
\A(h,\tilde{C}) \Delta(h,c) 
\simeq
\frac{\vol(s\K)}{\disc(\tilde{C}L)}
\sum_{h \mod \tilde{C}}
\A(h,\tilde{C}) \Delta(h,c),
$$
which is made rigorous by:
\begin{proposition}
If $\disc(L) \cdot \tilde{C} \leq s$ then
$$
\sum_{h \in s \K \cap L}
\A(h,\tilde{C}) \Delta(h,c) 
=
\frac{\vol(s\K)}{\disc(\tilde{C}L)}
\sum_{h \mod \tilde{C}}
\A(h,\tilde{C}) \Delta(h,c) 
+O(\tilde{C}  c^{-1/2+\epsilon}s^{r-2}  )
$$
\end{proposition}
\begin{proof}
  See 6.10 in \cite{part1}.
\end{proof}

Summing over $c,g$ and $L$ we get:
\begin{corollary}
$$
\frac{s}{ 2^{r \omega(q)}} 
\sum_{ \substack{ c|q \\ c \leq s^{1/3} \\ \omega(c) \leq \sqrt{\omega(q)} } }
\sum_{ \substack{ g | \frac{q}{c}  \\ g \leq s^{\Cbin} } }
\sum_{ \substack{ L \in \LL \\ \supp(L) = g \\ \disc(L) \leq s^{1/3}  } }
\lambda(L)
\sum_{h \in s \K \cap L}
\A(h,\tilde{C}) \Delta(h,c) 
$$
$$=
\frac{s}{ 2^{r \omega(q)}} 
\sum_{ \substack{ c|q \\ c \leq s^{1/3} \\ \omega(c) \leq \sqrt{\omega(q)} } }
\sum_{ \substack{ g | \frac{q}{c}  \\ g \leq s^{\Cbin} } }
\sum_{ \substack{ L \in \LL \\ \supp(L) = g \\ \disc(L) \leq s^{1/3}  } }
\lambda(L)
\frac{\vol(s\K)}{\disc(\tilde{C}L)}
\sum_{h \mod \tilde{C}}
\A(h,\tilde{C}) \Delta(h,c) 
$$
$$
+
O\left(
  \frac{s}{ 2^{r \omega(q)}} 
  \sum_{ \substack{ c|q \\ c \leq s^{1/3} \\ \omega(c) \leq \sqrt{\omega(q)} } }
  \sum_{ \substack{ g | \frac{q}{c}  \\ g \leq s^{\Cbin} } }
  \sum_{ \substack{ L \in \LL \\ \supp(L) = g \\ \disc(L) \leq s^{1/3}  } }
  |\lambda(L)| \tilde{C}  c^{-1/2+\epsilon} s^{r-2} 
\right).
$$
\end{corollary}
\begin{proof}
Immediate since the bounds on $c, \omega(c)$ and $\disc(L)$ forces
$\tilde{C} \cdot \disc(L)$ to be smaller than $s$ by
lemma~\ref{l:bound-d-tilde}. 

\end{proof}

\subsection{Estimating the error term}
By lemma~\ref{l:bound-d-tilde}, $\tilde{C}  c^{-1/2+\epsilon} \leq c
s^{1/6+\epsilon} \leq s^{1/2+\epsilon}$. This together with
lemma~\ref{l:lattice-bound} gives that the error term is bounded by 
$$
\frac{s^{r-1}}{ 2^{r \omega(q)}} 
  \sum_{ \substack{ c|q \\ c \leq s^{1/3} \\ \omega(c) \leq \sqrt{\omega(q)} } }
  \sum_{ \substack{ g | \frac{q}{c}  \\ g \leq s^{\Cbin} } }
  g^{2 \epsilon} s^{1/2+\epsilon} 
\ll
\frac{s^{r-1}}{ 2^{r \omega(q)}} 
  \sum_{ \substack{ c|q \\ c \leq s^{1/3} \\ \omega(c) \leq \sqrt{\omega(q)} } }
  \sum_{ \substack{ g | \frac{q}{c}  \\ g \leq s^{\Cbin} } }
  s^{1/2+2\epsilon} 
$$
which, by lemma~\ref{l:div-sum-bound}, is 
$$
\ll 
\frac{s^{r-1}}{ 2^{r \omega(q)}}
s^{1/2 +3\epsilon  }
\ll 
s^{-1/2+3\epsilon}
$$
and can thus be neglected.

\subsection{The main term}
In order to evaluate the main term we need to complete the sum,
i.e. extend it to all lattices and divisors:

\begin{lemma}
There exists $\CcompleteAgain > 0$ such that 
$$
\frac{s}{ 2^{r \omega(q)}} 
\sum_{ \substack{ c|q \\ c \leq s^{1/3} \\ \omega(c) \leq \sqrt{\omega(q)} } }
\sum_{ \substack{ g | \frac{q}{c}  \\ g \leq s^{\Cbin} } }
\sum_{ \substack{ L \in \LL \\ \supp(L) = g \\ \disc(L) \leq s^{1/3}  } }
\lambda(L)
\frac{\vol(s\K)}{\disc(\tilde{C}L)}
\sum_{h \mod \tilde{C}}
\A(h,\tilde{C}) \Delta(h,c) 
$$
$$
=
\frac{s}{ 2^{r \omega(q)}} 
\sum_{ c|q }
\sum_{ \substack{ g | \frac{q}{c} }}
\sum_{ \substack{ L \in \LL \\ \supp(L) = g }}
\lambda(L)
\frac{\vol(s\K)}{\disc(\tilde{C}L)}
\sum_{h \mod \tilde{C}}
\A(h,\tilde{C}) \Delta(h,c)$$
$$
+ O\left( \exp \left( -\CcompleteAgain \sqrt{\omega(q)}  \right) \right).
$$
\end{lemma}
\begin{proof}
By lemma~\ref{l:delta-sup} and corollary~\ref{c:eps-bound},
$$
\frac{1}{\tilde{C}^{r-1}} 
\sum_{h \mod \tilde{C}} 
\A(h,\tilde{C}) \Delta(h,c) \ll c^{1/2+\epsilon}
$$ 
and we can therefore use the same bounds as in lemma~\ref{l:chop-l} to
include the terms for which $\disc(L) \geq s^{1/3}$. Since
$$
\sum_{ \substack{ L \in \LL \\ \supp(L) = g  } } 
\frac{\lambda(L)}{\disc(\tilde{C}L)}
\ll g^{2 \epsilon-1}
$$ we can use lemma~\ref{l:div-sum-bound} to include $g \geq
s^{\Cbin}$. Finally, similar bounds used in lemma~\ref{l:chop-c}
allows us to extend the sum to include all $c,\tilde{C}$.
\end{proof}

The completed sum is multiplicative, and we can evaluate it as
follows: Expanding $N(h, \tilde{Q})$ we see that 
$$
\frac{1}{N_{\tilde{Q}}} 
\sum_{h \in (\Z/\tilde{Q}\Z)^{r-1}} N(h, \tilde{Q}) 
=
\frac{s_{\tilde{Q}}}{ 2^{r \omega(q)}} 
\sum_{ c|q }
\sum_{ \substack{ g | \frac{q}{c} }}
\sum_{ \substack{ L \in \LL \\ \supp(L) = g }}
\lambda(L)
\sum_{ \substack{ h \mod \tilde{Q} \\ h \in L  }}
\A(h,\tilde{C}) \Delta(h,c). 
$$
Since $\disc(L)$ and $c$ are coprime the intersection of a fundamental
domain of $\tilde{C}L$ with $L$ consists of a full set of
representatives of $\Z^{r-1}/\tilde{C}\Z^{r-1}$ (see lemma 8 in
\cite{part1}.) Now, $\R^{r-1}/\tilde{Q}\Z^{r-1}$ can be expressed as a
disjoint union of $\frac{\tilde{Q}^{r-1}}{\disc(\tilde{C}L)}$
translates of the fundamental domain for $\tilde{C}L$, and thus 
$$
\sum_{ \substack{ h \mod \tilde{Q} \\ h \in L  }}
\A(h,\tilde{C}) \Delta(h,c)
=
\frac{\tilde{Q}^{r-1}}{\disc(\tilde{C}L)}
\sum_{h \mod \tilde{C}}
\A(h,\tilde{C}) \Delta(h,c).
$$
Hence 
$$
\frac{1}{N_{\tilde{Q}}} 
\sum_{h \in (\Z/\tilde{Q}\Z)^{r-1}} N(h, \tilde{Q}) 
$$
$$
=
\frac{s_{\tilde{Q}}}{ 2^{r \omega(q)}}
\sum_{  c|q }
\sum_{ \substack{ g | \frac{q}{c} }}
\sum_{ \substack{ L \in \LL \\ \supp(L) = g }}
\lambda(L)
\frac{\tilde{Q}^{r-1}}{\disc(\tilde{C}L)}
\sum_{h \mod \tilde{C}}
\A(h,\tilde{C}) \Delta(h,c).
$$
On the other hand, $$\sum_{h \in (\Z/\tilde{Q}\Z)^{r-1}} N(h,\tilde{Q}) =
N_{\tilde{Q}}^{r}$$ since all $r$-tuples of squares are accounted for
when we sum over all $h$. Hence
$$
\frac{s_{\tilde{Q}}}{ 2^{r \omega(q)}}
\sum_{  c|q }
\sum_{ \substack{ g | \frac{q}{c} }}
\sum_{ \substack{ L \in \LL \\ \supp(L) = g }}
\frac{\lambda(L)}{\disc(\tilde{C}L)}
\sum_{h \mod \tilde{C}}
\A(h,\tilde{C}) \Delta(h,c)
=
\left(  \frac{N_{\tilde{Q}}}{\tilde{Q}} \right)^{r-1}
=
\frac{1}{s_{\tilde{Q}}^{r-1}},
$$
and thus
$$
\frac{s}{ 2^{r \omega(q)}} 
\sum_{  c|q } 
\sum_{  g | \frac{q}{c} }
\sum_{ \substack{ L \in \LL \\ \supp(L) = g }}
\lambda(L)
\frac{\vol(s\K)}{\disc(\tilde{C}L)}
\sum_{h \mod \tilde{C}}
\A(h,\tilde{C}) \Delta(h,c) 
$$
$$ =
\frac{ \vol(\K) s^r }{ 2^{r \omega(q)}} 
\sum_{  c|q } 
\sum_{  g | \frac{q}{c} }
\sum_{ \substack{ L \in \LL \\ \supp(L) = g }}
\frac{\lambda(L)}{\disc(\tilde{C}L)}
\sum_{h \mod \tilde{C}}
\A(h,\tilde{C}) \Delta(h,c) 
$$
$$
=
\vol(\K) \frac{s^r}{s_{\tilde{Q}}}
\frac{1}{s_{\tilde{Q}}^{r-1}} 
=
\vol(\K) 
\left( \frac{s}{s_{\tilde{Q}}} \right)^r,
$$ 
and we will be done if we can show that $\frac{s}{s_{\tilde{Q}}} = 1 +
O( \exp( -C \sqrt{\omega(q)})$ for some $C>0$. Now, it is easy to see that 
$$N_{p^k} = p^{k-1} \frac{p-1}{2} +
p^{k-1 -2} \frac{p-1}{2} + \ldots
+ \frac{p-1}{2} p^{k-1 - 2 (\lceil k/2 \rceil-1)} +1,
$$
which implies that 
$$
\frac{N_{p^k}}{p^k} = \frac{\sigma(p)}{2}
\left( 1 + p^{-2} + \ldots + p^{-2 (\lceil k/2 \rceil-1)}
  +\frac{2}{\sigma(p)} p^{-k}
\right).
$$
We may assume that $ \tilde{\alpha}_p < \alpha_p$, hence 
$$
\frac{s}{s_{\tilde{Q}}} = \prod_{p|q} \left(1+O(p^{-\tilde{\alpha}_p})
  \right)
=
\exp \left( \sum_{p|q} \log \left( 1+ O(p^{-\tilde{\alpha}_p}) )
  \right) \right)
=
\exp \left( \sum_{p|q} O(p^{-\tilde{\alpha}_p}) \right).
$$
But 
$$\sum_{p|q} p^{-\tilde{\alpha}_p} \ll 
\sum_{p|q} p^{-1/2} \exp \left( -\frac{ \sqrt{\omega(q)}}{\CexpOne} \right) \ll
\omega(q) \exp \left( -\frac{ \sqrt{\omega(q)}}{\CexpOne} \right)
$$
$$
\ll \exp \left( -\frac{\sqrt{\omega(q)}}{\CexpTwo} \right).
$$ 
Thus 
$$
\frac{s}{s_{\tilde{Q}}} = 
1 + O \left( \exp \left( -\frac{\sqrt{\omega(q)}}{\CexpTwo} \right) \right)
$$
and we have proved theorem~\ref{t:main}.
\appendix
\section*{Appendix}

\begin{lemma}
\label{l:f-sub-q-bound}
Let $p_1$ be the smallest prime dividing $q$. With $F(q,t) =
\sum_{p|q} p^{-t}$ and $k>0$ an integer we have 
\[
F(q,k/2) \leq 
\left\{
\begin{array}{ll} 
O \left( \sqrt{\frac{\omega(q)}{\log \omega(q)}} \right) & \mbox{ if $k=1$,} \\
O \left( \log( \log( \omega(q) )) \right) & \mbox{ if $k=2$, } \\
3 p_1^{1-k/2}  & \mbox{ if $k \geq 3$. }
\end{array}
\right.
\]
\end{lemma}
\begin{proof}
For the cases $k=1$ and $k=2$ we may assume that $q$ is the product of
the first $\omega(q)$ primes, and the bounds are then immediate
consequences of the prime 
number theorem, together with the fact that the $\omega(q)$-th prime
is roughly of 
size $\omega(q) \log \omega(q)$. For $k \geq 3$, we note that the sum is
bounded by $p_1^{-k/2} + \int_{p_1}^\infty x^{-k/2} dx < p_1^{-k/2} +
(k/2-1)^{-1} p_1^{1-k/2} < 3 p_1^{1-k/2}$. 
\end{proof}

\begin{corollary}
\label{c:eps-prod} There exist $\CepsProd >0$ such that 
$$
\prod_{p|q} (1 + \Ceps p^{-1/2}) 
\ll
\exp \left( \CepsProd \sqrt{  \frac{\omega(q)}{\log(\omega(q))} } \right).
$$
\end{corollary}
\begin{proof}
$$
\prod_{p|q} (1 + \Ceps p^{-1/2}) =
\exp \left( \sum_{p|q} \log (1 + \Ceps p^{-1/2}) \right) 
$$
$$\ll
\exp \left( \sum_{p|q} \Ceps p^{-1/2} \right) =
\exp \left( \Ceps F(q,1/2) \right).
$$
\end{proof}

\begin{lemma}
\label{l:power-series} There exists $C>0$ such that 
$$
\sum_{\substack{c|q  \\  \omega(c) \geq \sqrt{\omega(q)}  } }
\Ceps^{\omega(c)}
c^{-1/2} 
\ll_{\Ceps}
\exp( -C \sqrt{\omega(q)} ).
$$
\end{lemma}
\begin{proof} 
Let
$$
f(z)=
\prod_{p|q}
(1+z \Ceps p^{-1/2})
=
\sum_{k=0}^{\omega(q)} z^k a_k
$$
where $a_k = \sum_{c|q, \omega(c)=k} \Ceps^{\omega(q)} c^{-1/2}$.  By
Cauchy's theorem 
$$
a_n 
= 
\frac{1}{2 \pi i}
\int_{|z|=2}
\frac{ f(z) } {z^{n+1}} dz
$$
and thus
$$
|a_n|
\leq
\frac{1}{2 \pi }
\int_{|z|=2}
\frac{ |f(z)|} {2^{n+1}} |dz|.
$$

Write $q = q_1 \cdot q_2$
where $q_1 = \prod_{ \substack{  p|q \\ p \leq 9 \Ceps^2}  } p$,
$q_2 = \prod_{ \substack{  p|q \\ p > 9 \Ceps^2} } p$ and let 
$$f_i = \prod_{p|q_i} (1+z \Ceps p^{-1/2}).$$ Clearly $|f_1(z)|
\ll_{\Ceps} 1$ for $|z| \leq 2$. Moreover,
$$
f_2(z) =
\prod_{ p|q_2}
(1+ z \Ceps  p^{-1/2})
=
\exp \left(
\sum_{ p|q_2 }
\log(1+ z \Ceps  p^{-1/2})
\right)
$$
$$
=
\exp \left(
\sum_{k=1}^\infty 
\frac{(-1)^{k+1}}{k} 
F(q_2,k/2) 
( z \Ceps ) ^k 
\right).
$$
By lemma~\ref{l:f-sub-q-bound}, 
$F(q_2,k/2) \leq 3 (9 K_1^2)^{1-k/2}$ for $k \geq 3$ and thus $$
f_2(z)
=
\exp \left(
z \Ceps F(q_2,1/2) - \frac{(z \Ceps)^2}{2} F(q_2,1) + g(z)
\right)
$$
where $g(z)$ is an analytic function whose  $n$-th coefficient of
its power series expansion around zero is bounded by
$K_1^2 \cdot3^{3-n}$. Consequently, $|g(z)|$ is uniformly bounded in
$q$ as well as $z$ for $|z| \leq 2$. Hence 
$$
|a_n| \ll 
\frac{\exp \left( 2 \Ceps F(q,1/2) + \frac{(2\Ceps)^2}{2} F(q,1) \right)}
{2^{n+1}}
$$
for $|z| \leq 2$. By the bounds on $F(q,1/2)$ and $F(q,1)$ in
lemma~\ref{l:f-sub-q-bound} we see that 
$$
|a_n|
\ll
\frac{ \exp\left( O \left( \sqrt{ \frac{ \omega(q) }{ \log \omega(q)}}
    \right) \right)  } { 2^n}.
$$
Thus, 
$$
\sum_{\substack{ c|q \\  \omega(c) \geq \sqrt{\omega(q)} } }
\Ceps^{\omega(c)} c^{-1/2} 
=
\sum_{k = \sqrt{\omega(q)}}^{\omega(q)} a_k 
\ll
\exp \left( O \left( \sqrt{ \frac{ \omega(q) }{ \log \omega(q)}}
  \right) \right) \sum_{k =\sqrt{\omega(q)}}^\infty \frac{1}{2^{n}}
$$
$$
\ll
\exp \left( O \left( \sqrt{ \frac{ \omega(q) }{ \log \omega(q)}}
  \right) \right)
\frac{1}{2^{\sqrt{\omega(q)}}}
$$
$$
\ll
\exp \left( 
  O \left( \sqrt{ \frac{ \omega(q) }{ \log \omega(q)}} \right)
  -  \sqrt{\omega(q)} \log(2)  
\right)  \ll \exp(-C \sqrt{\omega(q)})
$$
for any $C < \log 2$. 
\end{proof}

\begin{lemma} \label{l:chop-complete} Let $f$ be a multiplicative
  function such that 
$f(c) \leq c^{-1/2} \Ceps^{\omega(c)}$ for some constant $\Ceps>0$. Then 
there exists constants $\CepsProd, \CchopC$ such that for all $q$
$$
\sum_{c|q} f(c)
=
\sum_{ \substack{ c|q \\ c \leq s^{1/3}  }  } f(c) + O(s^{-1/6+\epsilon})
=
\sum_{ \substack{ c|q \\ c \leq s^{1/3} \\ \omega(c) \leq
    \sqrt{\omega(q)} }  } f(c) +  
O\left( \exp \left( -\CchopC \sqrt{\omega(q)}  \right) \right) 
$$ 
and
$$
\sum_{c|q} f(c) \ll 
\exp \left( \CepsProd \sqrt{  \frac{\omega(q)}{\log(\omega(q))} } \right).
$$

\end{lemma}
\begin{proof}
For the first assertion we note that $f(c) \leq c^{-1/2}
\Ceps^{\omega(c)}$ implies that  
$$
\sum_{ \substack{ c|q \\ c \geq  s^{1/3}  }  } |f(c)|
\ll 
\sum_{ \substack{ c|q \\ c \geq  s^{1/3}  }  } 
c^{-1/2} \Ceps^{\omega(c)}
\ll
\sum_{ \substack{ c|q \\ c \geq  s^{1/3}  }  } 
c^{-1/2 +\epsilon},
$$
which by lemma~\ref{l:div-sum-bound} is bounded by
$s^{-1/6+\epsilon}$. The second assertion follows from
lemma~\ref{l:power-series}, and the last follows from
corollary~\ref{c:eps-prod}. 
\end{proof}

\begin{lemma}
\label{l:div-sum-bound}
Let $\alpha,\beta >0$. Then 
$$
\sum_{ \substack{ c|q \\ c \geq s^{\alpha}  }  } c^{-\beta}
\ll 
s^{-\alpha \beta + \epsilon }.
$$
Moreover, 
$$
\sum_{ \substack{ c|q \\ c \leq s^{\alpha}  }  } 1 \ll
s^{\epsilon}.
$$
\end{lemma}
\begin{proof}
See lemma~18 and 19 in \cite{part1}.
\end{proof}

{\em Acknowledgement.} The author wishes to thank Ze\'{e}v Rudnick for
helpful discussions and comments.

\end{document}